\documentclass[letterpaper, 10pt, conference]{ieeeconf}
\IEEEoverridecommandlockouts

\usepackage{graphicx}
\usepackage{float}
\usepackage{subfigure}
\usepackage{derivative}
\usepackage[utf8]{inputenc}
\usepackage{nicefrac}

\usepackage{bm}
 \newcommand{\manu}{manufacturer }
 \newcommand{\stack}{Stackelberg }
\usepackage{amsmath,mathptmx}

\usepackage{amssymb}
\usepackage{tcolorbox}
\usepackage{tikz}

\usepackage{mdframed}

\usepackage{dsfont}
\usepackage{graphicx}
\usepackage[colorlinks=true, allcolors=blue]{hyperref}

\newcommand{\D}{{  D}}

\newcommand{\bpi}{{\bm \pi}}

\newcommand{\br}{{\bm r}}
\newcommand{\R}{{\cal R}}
\newcommand{\Q}{{\cal Q}}

\newcommand{\bq}{{\bf q}}
\newcommand{\bp}{{\bf p}}
\newcommand{\bU}{{\bar {U}}}
\newcommand{\brid}{{\bf{r}}^{id}}

\newcommand{\Rev} [1]{{\color{black}  #1}}

\newcommand{\hPhi}{{\widehat \Phi}}
\newcommand{\eop}{{\hfill $\Box$}}
\newtheorem{thm}{{\bf Theorem}}
\newtheorem{lem}{{\bf Lemma}}

\newcommand{\ignore}[1]{}

\IEEEoverridecommandlockouts                              

\newcommand{\TR}[2]{#2} 

\usepackage{bbm}
\usepackage{enumerate}

\begin{document}
\title{\bf Punitive policies to combat misreporting in dynamic 
supply chains
}

\author{Madhu Dhiman, Atul Maurya,  Veeraruna Kavitha and Priyank Sinha \thanks{Department of Industrial Engineering and Operations Research, Indian Institute of Technology, Bombay. Emails: \{madhu.dhiman, 23m1517, vkavitha and priyank.sinha\}@iitb.ac.in}}
\maketitle
\begin{abstract}
Wholesale price contracts are known to be associated with double marginalization effects, which prevents supply chains from achieving their true market share. In a dynamic setting under information asymmetry, these inefficiencies manifest in the form of misreporting of the market potential by the manufacturer to the supplier, again leading to the loss of market share. We pose the dynamics of interaction between the supplier and manufacturer as the Stackelberg game and develop theoretical results for optimal punitive strategies that the supplier can implement to ensure that the manufacturer truthfully reveals the market potential in the single-stage setting. Later, we validate these results through the randomly generated, Monte-Carlo simulation based numerical examples.    
\end{abstract}

\ignore{
In supply chains with multiple agents having   asymmetric information, the agents with more information can resort to  misreporting,  to secure more favorable outcomes. This work investigates  one such scenario involving a single manufacturer and a supplier. Here the manufacturer  possesses accurate knowledge of the market potential and the supplier banks on this information to quote a dynamic price which can potentially lead to a more efficient supply chain. 
But the manufacturer  may have an incentive to misreport  and when the supplier sets  prices (for raw material)  based on this inaccurate information, the efficiency and the profit of the supplier can degrade substantially.

To counter this, we propose a punitive pricing strategy  for the supplier to discourage misreporting and ensure truthful reporting from the manufacturer. We model the interaction using a Stackelberg game and introduce a penalty term that the supplier applies to ensure no misreporting. The aim of the strategy is to ensure that the manufacturer maximizes its utility by truthfully reporting the market potential. Through rigorous analysis, we provide the punitive pricing strategy achieves optimal truth-revealing behavior.}

\vspace{-4mm}
\section{Introduction}
\label{Sec_intro}

In supply chains (SCs),   manufacturers buy raw material from the suppliers and  produce final products that can be sold to  end-customers. 
We consider a market with one manufacturer and one supplier, where the manufacturer can estimate the (randomly fluctuating) market potential accurately (\cite{vosooghidizaji2020supply}). 
It informs the supplier about the same to seek a dynamic price  for the raw material.  
The prices quoted by the manufacturer to the end-customers depend on various factors  (like cost of manufacturing, the 
market potential, etc.), including the price of the raw-material (see, e.g., \cite{wadhwa2024partition, zheng2021willingness, li2023pricing}); these prices dictate the actual demand attracted, which also drive the profits of the supplier. Thus, the supplier will be willing to  quote a dynamic price  (per unit of raw material) depending on the values of the market potential   reported by the manufacturer.

Due to asymmetry of information, the  \manu might have an incentive to misreport \Rev{(e.g., see \cite{zong2022decision})} to secure more favorable terms from the supplier; for example,  it is clear from the  results of \cite{wadhwa2024partition} that the supplier quotes a lower price when the market potential is small (this is true when the production and the procurement costs do not alter with  variations in the market potential, which in turn is true when the variations are not exuberant).   In such cases, the supplier is at the risk of being misinformed, potentially incurring losses by providing raw materials under sub-optimal dynamic prices.
\ignore{
\Rev{ There are many case studies to this effect, see   e.g.,  \cite{ravichandran2008managing}, \cite{anderson2000upstream}; in US during .. have reportedly inflated... \\
n....\\
There are many case studies to this effect, see   e.g.,  \cite{ravichandran2008managing}, \cite{anderson2000upstream}; in US during .. have reportedly inflated... 
}}
\Rev{ In many strands of literature and news articles, several instances of misreporting are discussed, (e.g., in  \cite{ravichandran2008managing} bullwhip effects are shown be resultant of misreporting, \cite{anderson2000upstream}. 
 Australia's largest building materials supplier, Boral,  identified cost misreporting by \$20–30 million, see \cite{Boral}, \cite{Reuters}, etc). 

}

 This work examines the intricate dynamics where the manufacturer reports the time varying market potentials to the supplier, who leads by dynamically setting the  prices (for raw-materials) based on the reported values. The main focus is   on the  design of the punitive policies for the supplier that compel the manufacturer to report truthfully. \Rev{We  discuss two types of punitive policies  and  theoretically establish their effectiveness.} We also provide an elaborate numerical study to illustrate that the policy implemented with intermediate estimates also succeeds in eliciting the truth from the  manufacturer. This work is inspired by the punitive policies proposed in \cite{kavitha2012opportunistic} in the context of wireless networks. 

 There are strands of literature that consider  asymmetric information games in SCs in other contexts (e.g., \cite{ni2021supply, zheng2021willingness,shen2019review}). 
For example, in \cite{ni2021supply}, the suppliers are not aware of the innovation efficiency (that drives the production costs) of the \manu and aim to design appropriate contracts. 
In a recent survey \cite{vosooghidizaji2020supply} on information asymmetry in SCs, it is observed that the majority of the literature   on demand information asymmetry focuses either   on designing optimal policies or on finding feasible solutions, while our work focuses on designing punitive policies that successfully elicit truthful information. By virtue of the proposed policies, the asymmetric information SCs can work efficiently as full information SCs.  We theoretically prove truth-elicitation   for a  greedy \manu who under-reports the potential to corner a lower price for raw materials. \Rev{We propose  another policy (which requires micro observations)  and prove that it   works successfully against a strategic \manu (over-reports   at some instances with smaller losses and under-reports at instances with much larger gains to confuse the supplier); we, in fact, prove theoretical results for the dynamic system.}

\newcommand{\RemEM}[1]{}  
\vspace{-2mm}
\section{Model Description}
\label{sec_model_desc}
\RemEM{ 
Consider a supply chain (SC) with one manufacturer and one supplier.
\RemEM{The customers purchase the final product from the manufacturer depending upon   factors like price, essentialness of the product, reputation of the manufacturer, etc. (see e.g., \cite{wadhwa2024partition}, \cite{zheng2021willingness}, \cite{li2023pricing}).} The \manu obtains the required raw materials from the supplier and sets a selling price for the final product and the demand depends upon this price  (see e.g., \cite{wadhwa2024partition}, \cite{zheng2021willingness}, \cite{li2023pricing}). The supplier, in turn, quotes a price for one   unit of  raw materials, required to produce one unit of final product. 
Both the agents  are rational and    aim   to maximize their individual revenue/shares by setting  appropriate prices, based on market conditions. The main focus of this paper is on the SCs facing  fluctuating market conditions. 
\RemEM{We begin with describing the key ingredients of the system, for fixed market condition.} 
}

In SCs, the manufacturer directly faces the customers and attempts to serve the maximum fraction of  them   at any given time point and at the best possible profit margin. 
It sets a selling price for the final product depending upon the market potential,  the price quoted by the supplier for the raw material and the production cost.  The supplier  quotes a price based on its own costs and the demand for the raw materials; this demand is  also dictated by  the market potential of the final product (see e.g., \cite{wadhwa2024partition,zheng2021willingness,  li2023pricing}). 

\noindent{\bf Demand-price  relation:} 
Let $\phi$ represent the potential of the market at any given time-point.  Technically, this represents 
the number of customers that demand for the product if the selling price is   zero.
We consider that the actual demand attracted  (out of $\phi$) by the \manu at selling price $p$ is given by the following linear relation (e.g.,    \cite{wadhwa2024partition}, \cite{zheng2021willingness}, \cite{li2023pricing}),
\begin{equation}
    D(p, \phi) = (\phi-\alpha p)^+, \mbox{ with } x^+:= \max\{0, x\}.
    \label{Eqn_manu_demand}
\end{equation}
In the above, $\alpha \in (0, 1]$ represents  the price-sensitivity parameter, 
$\alpha p$ is the fraction of the   customers lost by the manufacturer by setting the price at $p$. Also, 
observe
the demand   becomes zero  when    $p$ is higher than  $  \nicefrac{\phi}{\alpha}$.
 
\noindent{\bf Utilities:} 
Say the supplier quotes $q$ as the price for the  raw material and the manufacturer sets $p$ as the selling price. Then the manufacturer derives the following utility (see \eqref{Eqn_manu_demand}): 
\ignore{
To ensure the profitability and viability of operations within the supply chain, the market potential must exceed a defined threshold. The supplier can decide not to operate, or can
quote a price $q \in [0, \infty)$.
Specifically, the minimum market potential must satisfy the condition\footnote{The revised manufacturer’s utility incorporates an adjusted pricing strategy $\tilde{q}$, ensuring positivity for operability. As penalties increase $\tilde{q}$, the equilibrium price $p^*$ rises. The feasibility of this strategy is determined using the Stackelberg optimization framework, leading to the minimum market potential condition $d_i > \alpha(C_s + C_m)$.} $d_i > 2\alpha(\tilde{q} + C_m)  
\label{Dbarmin}$. This threshold ensures that the revenue generated from the demand attracted by the market potential is sufficient to cover both the operating and manufacturing costs. When the market potential fails to meet this threshold, it is no longer optimal for the agents to continue operations. In such cases, further analysis becomes redundant, as the agents would cease to operate due to insufficient profit margins. This concept is formally outlined in Lemma 4 of the Appendix (see \cite{wadhwa2024partition}), which establishes the importance of a minimum market potential as a critical determinant for the continuation of operations.
}

\vspace{-2mm}
{\small\begin{equation}
    U_M (p,q) := D(p, \phi)  (p-q-C_m)   
    = ( \phi-\alpha p)^+  (p-q-C_m), 
    \label{Eqn_manu_utility}
\end{equation}}%
where  $C_m$ is the production cost (per unit) incurred by the manufacturer. Basically,  the demand attracted by the manufacturer directly translates into the profit.  In fact, this demand also translates into the profit of the supplier: 

\vspace{-2mm}
{\small\begin{equation}
    U_S(p,q) := D(p, \phi) (q-C_s)  = (\phi-\alpha p)^+ (q-C_s),\label{Eqn_supplier_utility}
\end{equation}}%
here $C_s$ is the cost for procurement of a bundle of
raw material required for producing one unit of the final product.

Thus,  both the agents of  the SC derive higher utility by attracting a bigger (demand) fraction of  the potential~$\phi$. However, the individual revenues also depend upon the  two quoted  prices 
$p$ and $q$. This premise  leads to a non-cooperative game theoretic setting, where  cooperation between the two  agents plays a significant role. The key   element of cooperation stems from the desire to attract maximum demand out of  the  market potential (at any  given time-point). Further, it is well known that the market potential varies with time because of the market fluctuations (see e.g., \cite{ni2021supply}).
The main focus of this paper is to study the response of the two agents  to the fluctuating market. 
We begin with the background on the game for fixed market potential.


\ignore{
\subsection*{Suppliers' utility:} The price quoted by manufacturer to the customers depends on various factors such as  the cost of manufacturing the product, 
market conditions etc., along with the price of raw materials.  And the fraction of market potential realized or the fraction of customers buying the product depends significantly on the quote price,  along with other aspects like reputation(s) of the agent, and essentialness of the product etc.  This directly translates into the profits of   both the agents. This is the precise reason which motivates the supplier to quote a dynamic price depending upon the market potential indicated by the manufacturer. Thus, when supplier chooses to operate, then suppliers' utility is 
}

\noindent{\bf Stackelberg game for fixed market potential:} 
When the supplier first sets the price for the raw material and then the manufacturer responds to it, we have a Stackelberg game with the supplier as the leader and the manufacturer as the follower (as in \cite{wadhwa2024partition}, \cite{colson2007overview}, \cite{bard1998practical}). 
Such games for fixed market potential (for example) are solved in  \cite{wadhwa2024partition}   and we reproduce here the corresponding Stackelberg Equilibrium (SBE), for the purpose of completion.  This  game is  described by the utilities given in \eqref{Eqn_manu_utility} and \eqref{Eqn_supplier_utility}, and 
is solved for the case with perfect information (the agents have access to all the information),  under the following assumption:

\medskip
\noindent{\bf A.} The market potential is sufficiently large, 
$
\phi>\alpha (C_s + C_m). 
$

\medskip
 
By
\cite[Theorem 1]{wadhwa2024partition},  under assumption {\bf A},     there exists a unique 
SBE $(p^*, q^*)$ 
given by:

\vspace{-3mm}
{\small\begin{equation}
\label{Eqn_pstar_qstar_after_stackelberg}
p^*(\phi) = \frac{ \phi + \alpha(q^* (\phi) + C_m)}{2 \alpha}, \mbox{ and } q^* (\phi) = \frac{\phi + \alpha(C_s - C_m)}{2 \alpha}.
\end{equation}}
Further, the utilities at this SBE are given by

\vspace{-4mm}
{\small\begin{eqnarray}
  U_S^*  (p^*, q^*; \phi ) =2   U_M^* (  p^*, q^*; \phi ) = \frac{\left([ \phi-\alpha (C_s + C_m)]\right)^2}{8 \alpha}.
    \label{Eqn_manu_supplier_optimal_utility}
\end{eqnarray}}

From \eqref{Eqn_pstar_qstar_after_stackelberg}, the optimal prices $(p^*, q^*)$ are directly proportional to  the market potential $\phi$. As already mentioned, in reality,  the \manu faces the downstream market directly and hence can have perfect information\footnote{We assume perfect estimation and prediction  in this paper. A future study considering estimation/prediction errors would be more realistic. } about the true market potential $\phi$, while the supplier has to rely on the \manu for the same. Thus the \manu has incentive to misreport (under/over report) the market potential in a bid to procure raw materials at a lower price. \textit{Our main focus is to design an appropriate truth revealing market-potential-based price-policy for the supplier.}


\section{Quasi static market fluctuations } 
\label{subsec_quasi_static}
We now describe the fluctuating market.  We model the market fluctuations as quasi-static (such models are common in wireless context, e.g. \cite{kavitha2012opportunistic}) as explained below.

\begin{figure}
    \centering
    \vspace{-3mm}
    \includegraphics[width=0.6\linewidth, height=0.35\linewidth]{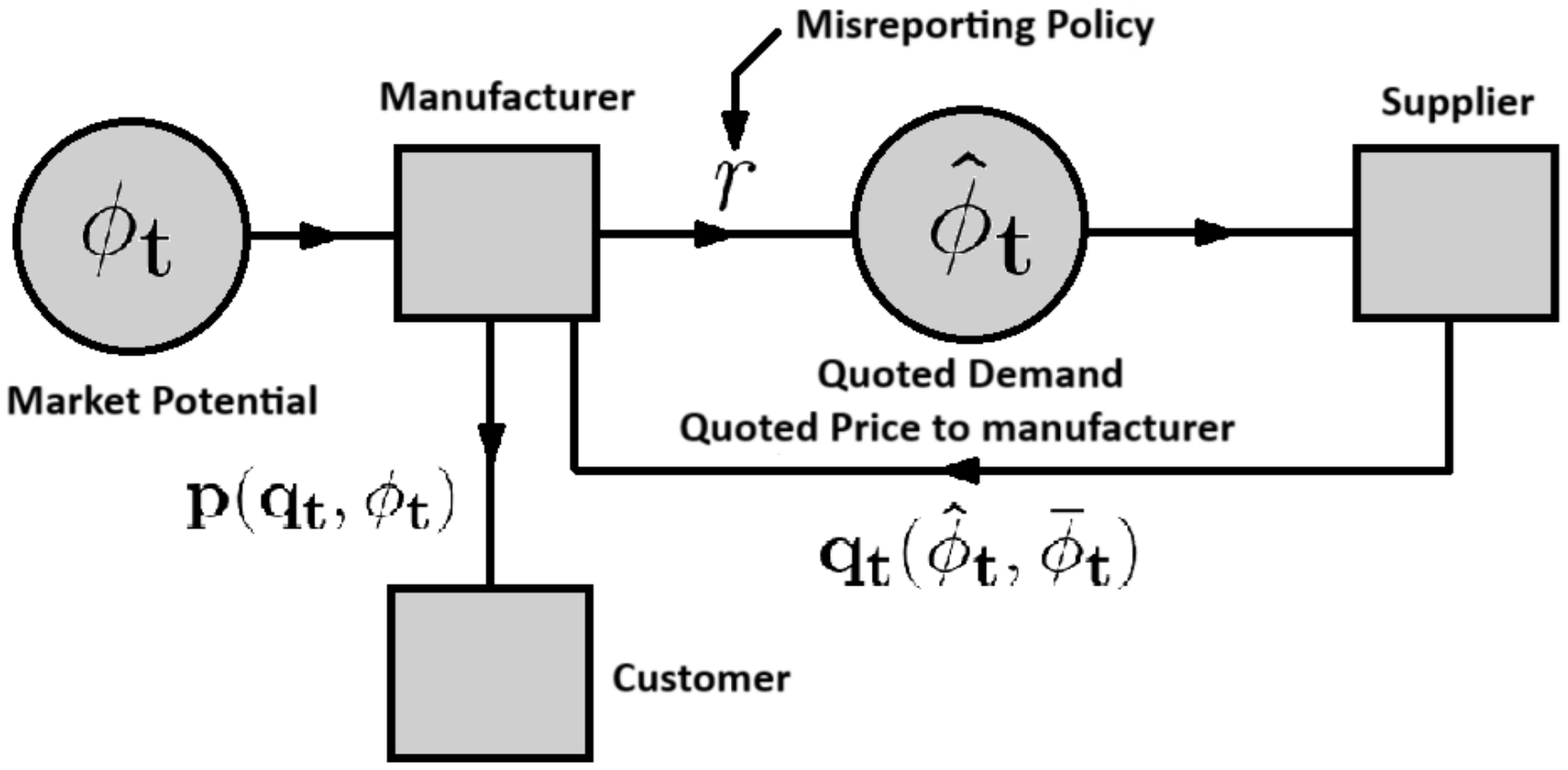}
    \vspace{-7mm}
    \caption{Supply Chain: flow of materials, prices and information}
    \label{fig:enter-label}
    \vspace{-5mm}
\end{figure}

The time frame is divided into slots of unit length (without loss of generality) and assume that the market potential is constant in each slot (each slot can represent a week, a month or any other time period specific to the market under consideration). Let $\Phi_t$ represent the potential in time slot~$t$, which can take   one of the  values from set $\Q = \{\phi_1, \phi_2, \cdots, \phi_L\}$ with respective probabilities $\{\sigma_1, \cdots, \sigma_L\}$; assume $\phi_l < \phi_{l+1}$ for all $l.$  Further,   $\{\Phi_t\}_t$   across different time slots are modeled to be  independent and identically distributed (IID). 

\subsection{Dynamic \stack game with perfect information}
\label{sec_limit_game}

The manufacturer faces the downstream market directly, we assume it to estimate the same accurately. For efficient operation of the SC, the supplier is willing to supply the raw materials at varying prices depending upon the potential of the slot. 
The following would be the flow of materials, information and prices/quotations in each slot. 
At the beginning of slot $t$, the manufacturer observes $\Phi_t$ and reports information about the same to the supplier; let $\hPhi_t$ represent this information. When the manufacturer reports   truthfully,  then $\hPhi_t  = \Phi_t$ (we assume this only in this subsection, \ref{sec_limit_game}). The supplier quotes price $Q_t = q(\hPhi_t)$   according to a Markov policy that we would describe in the immediate next. And, then the \manu sets the selling price at $P_t  =p(\Phi_t, Q_t)$ again according to  a Markov policy.

We consider stationary\footnote{It is well known that such policies are sufficient for sequential decision problems of this nature (see \cite{puterman1990markov}).} Markov (SM) policies, where the quoted price $q=q(\Phi)$ depends upon the (reported) market potential $\Phi$.  In other words, $\bq = (q_1, q_2, \cdots, q_L)$ represents a SM policy --- under this policy, the supplier quotes price $q_j$ in any time slot $t$, if the market potential in that slot  $\Phi_t =\phi_j.$
The policy as  
in~\eqref{Eqn_pstar_qstar_after_stackelberg} is  an example SM policy. 

The SM policy of the \manu   dictates the selling price of the final product given the market potential $\Phi_t$ and the price $Q_t$ quoted by the supplier. For SM supplier policies,  the policies of the \manu are also Markovian (under truthful reporting),  depending  only upon the potential $\Phi_t$  because of the following observation:   for any function $p (\cdot, \cdot)$ of the market potential and the supplier-price, we have: 

\vspace{-3mm}
{\small\begin{equation}
\label{Eqn_p_tildep}
  p(\Phi_t, Q_t) = p (\Phi_t, q(\Phi_t) ) 
= {\tilde p}(\Phi_t),
\end{equation}}%
for some appropriate function ${\tilde p} (\cdot)$  only of the market potential.
In other words, $\bp = (p_1, \cdots, p_L)$ represents a SM policy for the \manu ---  under this policy, the \manu sets the selling price of the product at $p_j$ in any time slot $t$ if the market potential in that slot equals $\Phi_t =\phi_j.$

In any time slot, the instantaneous utilities/revenues derived by both the agents are given by:

\vspace{-4mm}
{\small\begin{eqnarray}
   U_{M_t} &:=& U_{M} (P_t, Q_t; \Phi_t) = (\Phi_t - \alpha P_t)^+ (P_t - Q_t - C_m)  \mbox{ and } \nonumber \\
    U_{S_t} &:=& U_{S} (P_t, Q_t; \Phi_t) = (\Phi_t - \alpha P_t)^+ (Q_t - C_s).
    \label{Eqn_manu_supp_utility_slot_t}
\end{eqnarray}}
As the market potentials $\{\Phi_t\}_t$ across  
 various time slots are IID,    the pair of prices $\{ P_t, Q_t\}_t$ across the time slots are also IID under SM policies.  As a result, the sequences of instantaneous  utilities  $\{U_{M_t}\}_t$ and 
  $\{U_{S_t}\}_t$ are also IID. Hence, by Law of Large Numbers (LLNs), the limits of the time averages of  the  instantaneous revenues generated by the \manu  and the supplier almost surely  (a.s.) equal:

\vspace{-3mm}
{\small \begin{eqnarray}
    \frac{ \sum_{t \le T} U_M (P_t, Q_t; \Phi_t ) }{T} \hspace{-1.3mm} &\hspace{-1.3mm} \to\hspace{-1.3mm} &\hspace{-1.3mm}   E[U_M(p(\Phi), q(\Phi)); \Phi] \nonumber 
 = \hspace{-1.3mm}\sum_i \sigma_i U_M(p_i, q_i; \phi_i),   \nonumber \\
 \frac{  \sum_{t \le T} U_S (P_t, Q_t; \Phi_t )}{T}   \hspace{-1.3mm} &\hspace{-1.3mm} \to\hspace{-1.3mm} &\hspace{-1.5mm}   \sum_i \sigma_i U_S(p_i, q_i; \phi_i), 
 \label{Eqn_us_um_llm}
\end{eqnarray}}%
where $U_M(\cdot, \cdot)$  and $U_S(\cdot, \cdot)$ are defined in \eqref{Eqn_manu_supp_utility_slot_t}. 

\noindent{\textit {\bf Limit game:}}
  We now define a limit Stackelberg game, under SM policies, to study  the interactions between the supplier and the \manu (with supplier setting the prices first based on the instantaneous market potential) for the scenarios with   fluctuating market. When $T$,  the time horizon is sufficiently large, the utility of both the agents can be approximated using LLNs as  in \eqref{Eqn_us_um_llm}. Thus,  the limit  \stack game is defined by the following utility functions

  \vspace{-3mm}
{\small
\begin{equation}
    \bU_M (\bp, \bq) :=\sum_i \sigma_i U_M(p_i, q_i; \phi_i),  \mbox{ and, }    \
    \bU_S (\bp, \bq) :=  \sum_i \sigma_i U_S(p_i, q_i; \phi_i), 
    \label{Eqn_manu_supp_utility_bar}
\end{equation}}
with  the decision variables as the SM policies $\bp$ and $\bq$.

\ignore{
The following would be the flow of materials, information and prices/quotations in each slot. 
At the beginning of slot $t$, manufacturer observes $\Phi_t$ and reports information about the same to the supplier; let $\hPhi_t$ represent this information. Say $\hPhi_t$ provides perfect information about $\Phi_t$, i.e.,  say $\hPhi_t  = \Phi_t$. Then the supplier quotes price $Q_t = q^*(\hPhi_t)$   given in \eqref{Eqn_pstar_qstar_after_stackelberg}. And then the \manu sets the selling price at $P_t  =p^* (\Phi_t)$.  
Thus, if the information was transferred perfectly and if the supplier responds optimally, the time average of  the revenue generated by the \manu  equals almost surely  (by Law of Large numbers)
\begin{eqnarray*}
    \frac{1}{T}   \sum_{t \le T} U_M (P_t, Q_t; \Phi_t ) &\stackrel{T\to \infty}{\longrightarrow} &E[U_M(p^*(\Phi), q^*(\Phi))] \\
 &= &\sum_i \sigma_i U_S(p_i, q_i; \phi_i)  .   
\end{eqnarray*}
where $U_M(\cdot, \cdot)$ is defined in \eqref{Eqn_manu_supp_utility_slot_t} (throughout the paper, some symbols are skipped to improve the narrative when there is clarity).
The time average of the revenue generated by the supplier, under the same conditions, almost surely equals:
\begin{eqnarray*}
    \frac{1}{T}   \sum_{t \le T} U_S (P_t, Q_t; \Phi_t ) \hspace{1mm}\stackrel{T\to \infty}{\longrightarrow}  \sum_i \sigma_i \frac{\left([ \phi_i-\alpha (C_s + C_m)]\right)^2}{8 \alpha}.   
\end{eqnarray*}
}

For time varying market conditions, the assumption {\bf A} has to be satisfied with minimum potential $\phi_1$, i.e., 

 \medskip
\noindent{\bf A$'$.} We require, $\phi_1 > \alpha (C_s + C_m)$. 
 \medskip

Now consider a \stack (dynamic) game,  with utilities given by \eqref{Eqn_manu_supp_utility_bar}. 
By separability of component-wise functions in \eqref{Eqn_manu_supp_utility_bar} (by optimizing for each $p_i$ separately i.e.,  for each $i$ in the best responses and  for each $q_i$ separately in the top-level optimization problem) and proceeding as in \cite{wadhwa2024partition}, one can immediately prove the following result:
\begin{lem}{\bf [SBE for dynamic market]} Under {\bf A$'$}, the SBE and  the utilities and  is given by,

\vspace{-3mm}
  {\small  \begin{eqnarray}
p^*_i = \frac{ \phi_i + \alpha(q^*_i(\phi_i) + C_m)}{2 \alpha}, \mbox{ and } q^*_i = \frac{\phi_i + \alpha(C_s - C_m)}{2 \alpha} \ \forall i, 
\label{Eqn_pstar_qstar_limit_SBE}
 \\
   \bU_S^*  (\bp^*, \bq^* ) = 2   \bU_M^* (  \bp^*, \bq^* ) = \sum_i \sigma_i \frac{\left([ \phi_i-\alpha (C_s + C_m)]\right)^2}{8 \alpha} \label{Eqn_Um_star_limit_SBE} .  
  \hspace{2mm}  \mbox{ \eop}
\end{eqnarray}}  
\end{lem}
 \noindent When the supplier has no access to the market information, as already discussed, the \manu   may have incentive to misreport. The supplier may not be able to sense misreporting instantaneously, but 
 can sense  the  same  in long run. This forms the basis for designing a counter policy for the supplier (different  from  \eqref{Eqn_pstar_qstar_limit_SBE})  to combat misreporting.

\section{Misreporting and a Punitive Policy}
\label{sec_greedy_manu}
From \eqref{Eqn_pstar_qstar_limit_SBE}, the greedy \manu can misreport a lower market potential to obtain the raw materials at  a lower price. Let   \vspace{-3mm}
\begin{equation}  
    r_{ij} := \mathbb{P}(\hPhi = \phi_j \mid \Phi = \phi_i), \quad \forall \phi_i, \phi_j \in \Q, 
    \label{Eqn_rij}
\end{equation}  
  represent the probability of misreporting the market potential as $\phi_j$ when the actual  value  is $\phi_i$, with $i \neq j$. \textit{Observe, $ r_{ii} $ represents the probability of truthful reporting when the actual market potential is $ \phi_i $.}
Let the misreporting policy of greedy \manu be represented by  $\br\in \R_e$, where,

 \vspace{-4mm}
 {\small \begin{equation}
     \R_e := \{ \br = \{r_{ij}\}_{i,j}:    r_{ij} \ge  0,  \forall (i,j), \  r_{ij} = 0  \forall i < j,     \mbox{ and }  \sum_j r_{ij} = 1 \ \forall i   \}. \label{Eqn_Re}
\end{equation}}%
\textit{We  begin with 
a greedy manufacturer, that only under reports, by setting  $r_{ij} = 0$ for all $i < j$, 
to (always)   derive a better quoted price $q$ (see \eqref{Eqn_pstar_qstar_after_stackelberg}).}

The supplier is not aware of  the instantaneous market potentials $\{\Phi_t\}$, however it   knows the  statistics described by $\Q$  and $\{\sigma_l\}_l$. Thus, it can  compute the expected value of the market potential (using $\sum_l \sigma_l \phi_l$),   can estimate the expected value  of the reported potentials (using $\frac{1}{T} \sum_{t\le T} \hPhi_t$)  and can sense the (extent of) misreporting using the difference between the two. \textit{This is the key idea:  one can use the difference  to design an appropriate punitive policy for supplier}. 

\noindent{\bf Greedy Punitive Policy:} 
We once again consider the limit game   as in sub-section \ref{sec_limit_game} to design and analyze the  punitive policy. 
The numerical performance of the proposed   policy  in the real-time, is studied using   Monte-Carlo   simulations in  section \ref{sec_numerical}. \Rev{The dynamic system  is directly analyzed  in section \ref{sec_smart_manu}, while studying the second punitive policy.}

If the \manu uses misreporting policy $\br \in \R_e$, then the difference  equals (using \eqref{Eqn_rij} and \eqref{Eqn_Re}):

\vspace{-4mm}
{\small\begin{eqnarray}
    E[\Phi] - \frac{1}{T} \sum_{t\le T} \hPhi_t \to   f(\br) : =  \sum_{ (i, j) : i \ge j}  \sigma_i (\phi_i - \phi_j) r_{ij}, \mbox{ a.s.}
    \label{eqn_fr}
\end{eqnarray}}%
Here, $ f(\br)$ denotes the expected difference between the actual market potential  and the quoted value after long run, when \manu uses $\br$ as the misreporting policy.

The supplier utilizes the above difference term \eqref{eqn_fr}, to design the following punitive price policy:

\vspace{-3mm}
{\small\begin{eqnarray}
    {q}_{\pi}(\hat{\phi}) &=& q^*(\hat 
     \phi) + \frac{\pi}{\sqrt{\alpha}} f(\br), \mbox{ \normalsize where, }
    \label{Eqn_q_tilde}
\end{eqnarray}}%
$q^* (\cdot)$ is as defined in \eqref{Eqn_pstar_qstar_limit_SBE}; division by $\sqrt{\alpha}$ is for the purpose of using some simplified constants  in some latter computations.
\textit{Here $\pi$ is the parameter of the supplier policy.}  
We would like to stress here that in pre-limit, the supplier would  neither know the instantaneous deviations, nor would it know $\br$. However as $t$ increases, the supplier gets a good estimate of $f(\br)$ by LLNs. The numerical study of section \ref{sec_numerical} obtains the performance of the punitive policy implemented  using the estimates, 

\vspace{-5mm}
{\small 
\begin{eqnarray} \hspace{15mm}
    f_t :=   E[\Phi] - \frac{1}{t} \sum_{s\le t} \hPhi_s  \mbox{ in place of }f(\br). \label{Eqn_inst_policies}
\end{eqnarray}}%
Such a numerical study also establishes the effectiveness of the policy, even after considering the estimation  and convergence errors  \Rev{(Theorem \ref{Thm_policy_smart_manu} provides the theoretical confidence).}

\newcommand{\Uij}{{\mathbb C}}

In every time slot $t$, in response to the  price   $Q_t = q_\pi (\hPhi_t)$ quoted  by the supplier using policy $\pi$ given in \eqref{Eqn_q_tilde},  the \manu derives the optimal selling 
price  by solving:

\vspace{-2.5mm}
{\small\begin{equation}
    \arg \max_p  \left ( \Phi_t - p \alpha \right )^+ (p - Q_t - C_m). \label{Eqn_argmax_demand}
\end{equation}}%
The above function can be optimized using simple derivative based techniques and boundary conditions, and 
at the limit, the corresponding   conditional slot-wise optimal utilities   are: 

\vspace{-4mm}
{\small\begin{eqnarray}
    \Uij(i,j) \hspace{-1mm}&\hspace{-2mm}:=\hspace{-2mm}&\hspace{-2mm}   \sup_{p \in [0, \infty)} E[ U_m |\Phi_t = \phi_i,  \ {\hPhi}_t = \phi_j, \ Q_t = q= q_\pi(\phi_j) ]
   \nonumber  \\
    && \hspace{-15mm} =  
    \left(\left [ h_{ij} - \pi f(\br) \right ]^+\right)^2 , \mbox{ } \label{Eqn_u_ij} 
    h_{ij} \ := \ \frac{2\phi_i - \phi_j - \alpha (C_m+C_s)}{4 \sqrt{\alpha}}.  \hspace{2mm}\label{Eqn_h_ij}
\end{eqnarray}} 

\newcommand{\tU}{ {{\tilde U}^{P_I}}}

Thus, the utilities $\tU_{M}(\cdot, \cdot)$ \textit{of the limit game   with partial information}, when \manu  misreports using policy $\br$ 
and when the supplier uses punitive policy $\pi$ as in \eqref{Eqn_q_tilde}  are given by (with conditional utilities, $\Uij(i,j)$ as defined in \eqref{Eqn_u_ij}):

\vspace{-4mm}
{\small \begin{eqnarray}
     \tU_{M}(\br, \pi) &=&  \sum_{{i\ge j}} \sigma_{i} r_{ij} \Uij(i, j). \label{eqn_um_pi_rpi}
 \end{eqnarray}}
 Let $\brid \in \R_e$ represent the truthful policy of the manufacturer, i.e., without misreporting and hence with $r^{id}_{ii} = 1$, for all~$i$. One important and  interesting observation here is that $\tU_M(\brid, \pi) = \bU^*_{M}(\bp^*, \bq^*)$,  the optimal utility \eqref{Eqn_Um_star_limit_SBE} derived by the \manu in  perfect information game given in  \eqref{Eqn_manu_supp_utility_bar}; more interestingly this is   true for all  the punitive policies, irrespective of $\pi$,  and this is derived when  the \manu does not misreport (that is, uses~$\brid$).

Hence, the \manu is guaranteed to derive  utility as in  perfect information game
by using truthful policy $\brid$; in such a case,   the supplier also derives as in \eqref{Eqn_Um_star_limit_SBE}. 

 We call a punitive policy $\pi$ to be \textit{truth-revealing policy}, if 
\begin{eqnarray*}
    \max_{\br } \tU_M(\br, \pi) = \tU_M(\brid, \pi), \mbox{ and if }  
      \tU_M(\br, \pi) <\tU_M(\brid, \pi), 
\end{eqnarray*}
for all $ \br \ne \brid.$
 We will now prove the first  main result of the paper by showing the existence of truth-revealing policies among the punitive policies given in \eqref{Eqn_q_tilde}, against  all the greedy policies $\br \in \R_e$
 (proof is in Appendix):
 \begin{thm}
 \label{Thm_greedy_policy}
     There exits a $\bar{\pi}>0$ such that for all $\pi > \bar \pi$
     $$
     \tU_{M}(\br, \pi) < \tU_{M}(\brid, \pi), \mbox{ for all } \br  \ne \brid , \br \in \R_e . \hspace{8mm} \mbox{\eop}$$ 
 \end{thm}
\noindent Thus there exists a policy $ {\pi}$ such that the manufacturer will    find it strictly beneficial to reveal the true market potential to the supplier -- 
   policy \eqref{Eqn_q_tilde}  with $\pi > \bar \pi$ is effective in ensuring  truth revelation by the greedy manufacturer.

\vspace{-2mm}
\section{Numerical study}
\label{sec_numerical}

The punitive policies \eqref{Eqn_q_tilde} (with large enough $\pi$) are established to be truth-revealing for limit game by Theorem~\ref{Thm_greedy_policy}.  We now examine their performance for finite time horizon, by using finite-horizon estimates \eqref{Eqn_inst_policies} in place of  the limit value $f(\br)$ of \eqref{Eqn_q_tilde}. 
Towards this, 
we will implement the policies  of the manufacturer ($\br$) and  the supplier ($\pi$) using Monte-Carlo (MC) simulations as described next.   
An IID sequence of market potentials $\{\Phi_t\}_{t\le T}$ are generated using  the given distribution  $ \{\sigma_i\}$
 and $ {\cal Q} $. Generate the sequence   $\{{\hPhi}_t\}_{t\le T}$ using $\br$ and the sequence $\Phi_t$, which forms the  (IID) sequence of  reported potentials. The    time averages of the reported potentials are generated    iteratively using, 

\vspace{-2mm}
{\small$$
{\overline \Phi}_{t+1} = {\overline \Phi}_t + \frac{1}{t+1} \left (  {\hPhi}_{t+1}  - {\overline \Phi}_t \right ),
$$}
and then the  supplier-prices  are generated using: 

\vspace{-3mm}
{\small\begin{eqnarray}
Q_t =   q_t(\hPhi_t) 
  = q^*({\hPhi}_t) + \pi f_t \mbox{ with } f_t:=   \bigg (  \  \sum_{i} \sigma_i \phi_i - {\overline \Phi}_t  \bigg )^+, 
  \label{Eqn_q_t_numerical}
\end{eqnarray}}%
instead of~\eqref{Eqn_q_tilde}. Finally,  
the time averages of the revenues of the manufacturer is obtained using:

\vspace{-5mm}
{\small\begin{eqnarray}
    {\overline U}_{M,t+1} = {\overline U}_{M,t}+ \frac{1}{t+1} \left (   \left (  \left [  h_{\Phi_t, \widehat{\Phi}_t}  - \pi f_t \right ]^+ \right )^2 - {\overline U}_{M,t} \right ) .
    \label{Eqn_bar_u_numerical}
\end{eqnarray}}

\begin{figure}[htbp]
\vspace{-4mm}
  \hspace{-4mm}
    \centering
    \begin{minipage}{0.225\textwidth}
        \centering
                \includegraphics[trim = {2.1cm 6cm 1.2cm 2.cm}, clip, width = 4.1cm, height = 3.1cm]{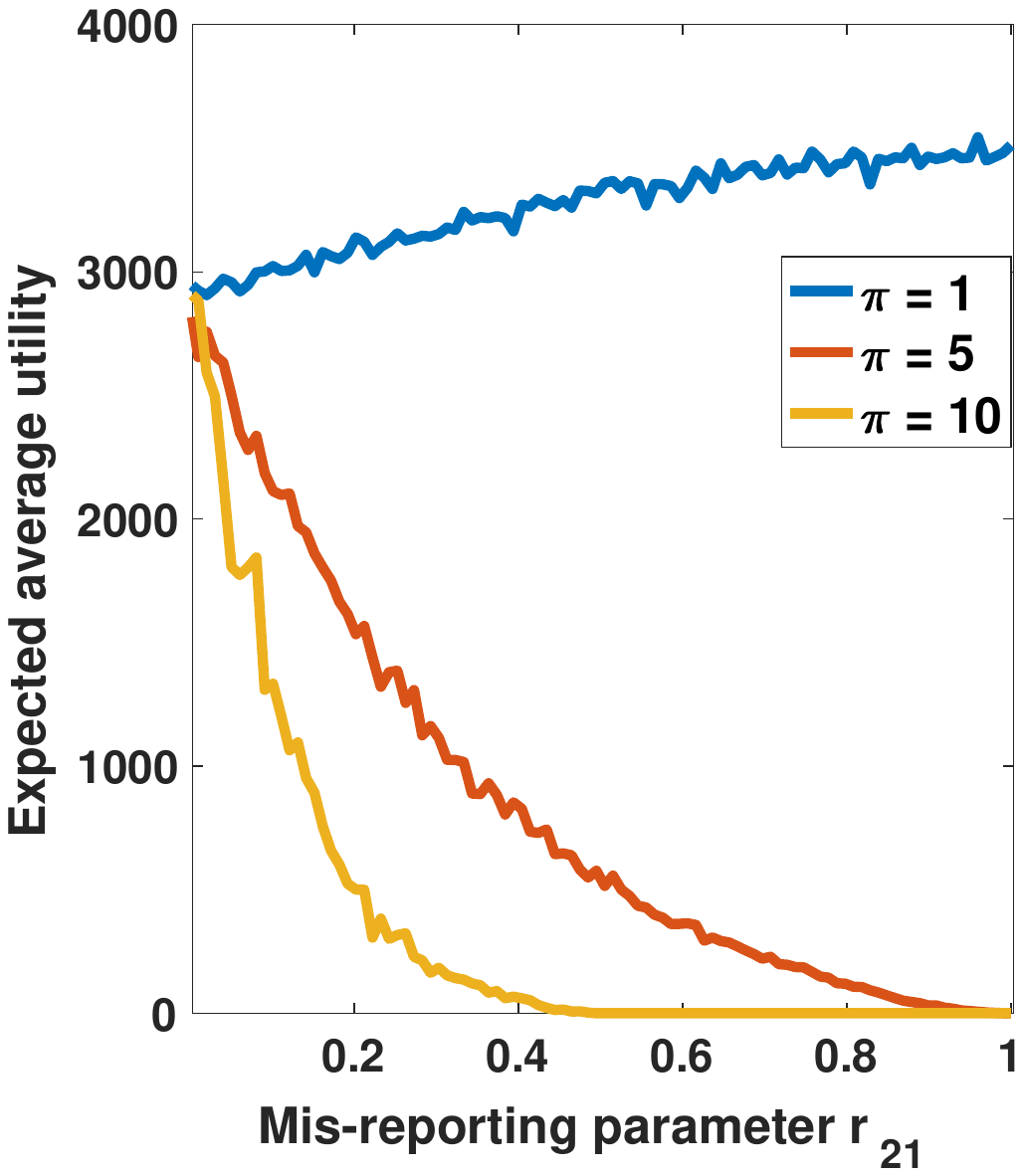}
    \end{minipage}
    \hspace{-2mm}
    \begin{minipage}
    {0.225\textwidth}
    \vspace{-2.5mm}
        \centering
       
            \includegraphics[trim = {2.1cm 6cm 1.8cm 3.cm}, clip, width = 4.4cm, height = 3.2cm]{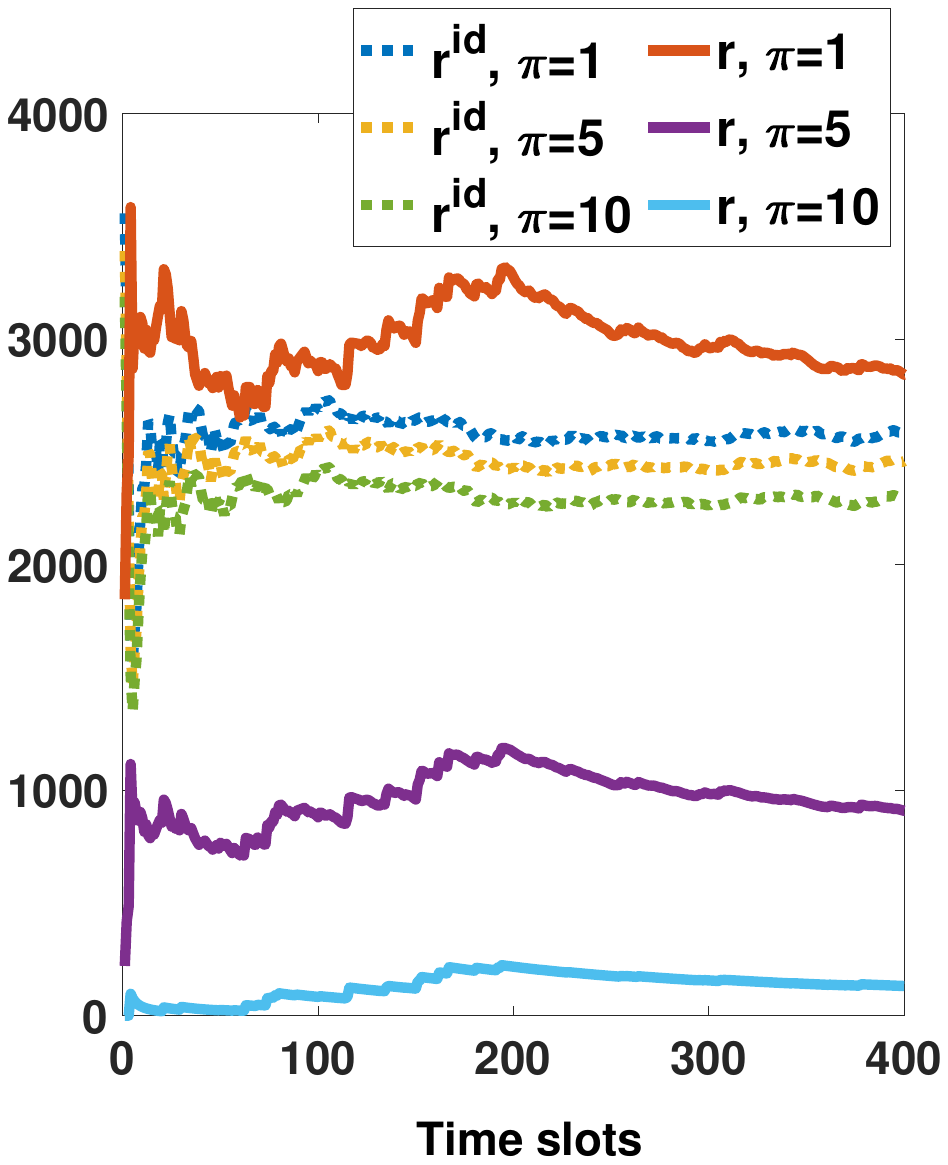}
        
    \end{minipage}
\vspace{-2mm}
    \caption{
    {\bf Left:} Impact of punitive policies on various misreporting policies 
    {\bf Right:} Time averaged utilities as a function of $T$ for different   $\br$ and $\pi$.   \label{fig:Manufacturer-Time-Averaged-Revenue-vs-r}
   }
\end{figure}


We consider a first example in the left sub-figure  of Figure~\ref{fig:Manufacturer-Time-Averaged-Revenue-vs-r} with   $C_s = 10$,   $C_m = 10$ and   $\alpha = 0.5$.  The number of time slots $T = 10000$,  the set  $\Q = \{40, 60, 80\}$ with probabilities $\sigma = \{0.2, 0.5, 0.3\}$. In this example,  we examine the impact of the misreporting policy on the time average utilities of the manufacturer; towards this, we set  $r_{32}=0$ and   fix  \(  r_{31} = r_{21}  \)  at a common value   and plot the utilities as a function of $r_{21}$; we plot three different curves   for three different values of $\pi$ (more details are provided in  the figure). 

\ignore{
    For (a), the analysis is conducted for different misreporting parameters $\pi \in \{1, 5, 10\}$. The market potential states are $\Phi = \{40, 60, 80\}$ with probabilities $\sigma = \{0.2, 0.5, 0.3\}$, and the misreporting policy is defined by $r_{21} = r_{31}$ and $r_{32} = 0$.  
We now examine the impact of the misreporting policy parameter \( r_{21} = r_{31} \) on the manufacturer's time-averaged revenue under varying penalty terms \( \pi \). The simulation results are presented in Figure \ref{fig:Manufacturer-Time-Averaged-Revenue-vs-r}, where the time-averaged revenue is plotted against the misreporting probability for different values of \( \pi \).  }

From the figure, it is evident that the time-averaged revenue increases with the misreporting probability $r_{21}$  at  lower value of penalty \(\pi = 1\) (see the blue curve). 
The upward trend in revenue   with  low penalty indicates that the supplier's punitive policy is not sufficiently strong to deter the manufacturer from engaging in misreporting.  
As the penalty term increases (\(\pi = 5, 10\)), the manufacturer's time-averaged revenue (red and orange curves) exhibits a clear downward trend with increasing misreporting probability $r_{21}$. This suggests that higher penalties effectively discourage  misreporting (and this is illustrated even after using estimates $f_t$ in place of the limit $f(\br)$). In fact the best utility for manufacturer is at the truth-eliciting $\brid$ (or at $r_{21}=0$).


For the second example of the right sub-figure of Figure~\ref{fig:Manufacturer-Time-Averaged-Revenue-vs-r}, 
we consider a market with larger variations, by setting   $\Q = \{10, 40, 60, 70, 80\}$ and probabilities $\sigma = \{0.1, 0.2, 0.3, 0.2, 0.2\}$; we have $T=400$ and the rest of the parameters are as before. 
We now analyze the manufacturer's time-averaged revenue under truthful and misreporting strategies with different penalty terms \( \pi \).  
Again at low penalty with \(\pi = 1\), given by red continuous curve,  the manufacturer achieves higher revenue by misreporting (with $\br \ne \brid$) than  that obtained under truthful policy $\brid$ (dotted curves). 
%
As the penalty increases (\(\pi = 5, 10\)), the revenue from truthful reporting $\brid$  exceeds that from misreporting, highlighting the effectiveness of stronger punitive policies. 



\section{Strategic Manufacturer and convergence}
\label{sec_smart_manu}
In  section \ref{sec_greedy_manu}, we have established that  policy   \eqref{Eqn_q_tilde} is effective against the greedy manufacturer. 
But such a policy will not work against a completely strategic manufacturer. Say the \manu over-reports for some market potentials (i.e., if $r_{ij} > 0$ for some $i < j$), and under-reports at some other potentials such that the expected deviation $f(\br) = 0$; then no value of $\pi$ can ensure truthful revelation.  The \manu can still benefit by quadratic nature of the instantaneous utilities at different market potentials (see \eqref{Eqn_u_ij}). 
In this case, the supplier has to 
design a more sophisticated policy, which probably involves  more micro observations (than just the deviations of the average of the  reported   and the actual   market potentials). \Rev{We now propose 
one such policy and analyze its  performance directly  for the dynamic system of section \ref{sec_numerical}; we in fact establish  its robustness theoretically.}

The supplier can observe the actual demand attracted in every time slot -- the demand attracted by the \manu in the entire duration of the time slot $t$,   directly translates to a proportional demand of the raw materials in the same slot. 
From \eqref{Eqn_manu_demand},  the demand attracted in time slot $t$   is given by 

\vspace{-3mm}
{\small \begin{equation}
    D_t  := D(P_t, \Phi_t) = \bigg ( \Phi_t - \alpha P_t(\Phi_t, Q_t(\hPhi_t))  \bigg )^+ + {\cal N}_t, 
    \label{Eqn_dt_smart_manu}
\end{equation}}%
\Rev{where ${\cal N}_t$ is some bounded  and zero mean IID  sequence that  captures the uncertainty in the demands realized
(here, $P_t (\phi, q) =  \nicefrac{(\phi+\alpha(q+C_m))}{2\alpha}$ represents the optimal choice of the manufacturer when the market potential  is $\phi$ and  the supplier price is $q$, essentially solves \eqref{Eqn_argmax_demand}.}

In other words, the supplier can estimate separately, for each $\phi_j$, the time average of the demands attracted in  the subset of the time slots when the reported value equals $ \phi_j$:

\vspace{-2mm}
{\small\begin{equation}
   {\bar D}_{j,t}   := \frac{ \sum_{s \le t} D_s \mathds{1}_{ \{ \hPhi_s = \phi_j \}}   } {  \sum_{s \le t} \mathds{1}_{ \{\hPhi_s = \phi_j \}}  }. \label{Eqn_report_demand_smart_manu} 
\end{equation}}%
Let ${\bar {\bf D}}_t := \{{\bar D}_{j,t}  \}_{ j }$ represent the vector of the above estimates. 
 If the \manu reports truthfully  by using $\brid$, the   quantities in \eqref{Eqn_dt_smart_manu}-\eqref{Eqn_report_demand_smart_manu}, converge     again by LLNs,  when\footnote{For SM   $\br$, $\{\hPhi_t\}$ is an IID sequence, we analyze only those  $j$ with $P(\hPhi = \phi_j) > 0$ as others are never reported.} $P(\hPhi_t = \phi_j) > 0$, see \eqref{Eqn_pstar_qstar_after_stackelberg} and \eqref{Eqn_h_ij}:

\vspace{-3mm}
{\small\begin{eqnarray*}
     {\bar D}^{id}_{j, t}  \stackrel{\mbox{as  $t\to \infty$}}{\rightarrow} 
\frac{\phi_j - \alpha(C_s+C_m)}{4}  =  \sqrt{\alpha } h_{jj}, \mbox{ \normalsize (a.s.)}. 
\end{eqnarray*}}%
 \Rev{The supplier can compute the above  and  any deviation from  $\sqrt{\alpha } h_{jj}$  indicates mis-reporting. We now propose a second punitive policy      using deviations $\{ {\bar { D}}_{j, (t-1)}  -  \sqrt{\alpha } h_{jj} \}_j$:   
 quote the following  price 
in the   $t$-th time slot:

 \vspace{-4mm}
{\small
\begin{eqnarray}
\label{Eqn_q_phi_smart}
Q_t (\hPhi_t, {\bf \bar D}_{t-1} )   = q^*(\phi_j) + 2{\pi}\left(  {\bar { D}}_{j, (t-1)}       -  \sqrt{\alpha }h_{jj} \right), \mbox{  \normalsize if }\hPhi_t = \phi_j. 
\end{eqnarray}}%
The above can again be seen as an SM policy by expanding the state  to include ${\bf \bar D}_{t-1}$. However,  one cannot derive asymptotic analysis of \eqref{Eqn_report_demand_smart_manu} just using  LLNs, when  $Q_t$ is as in \eqref{Eqn_q_phi_smart}.  Nonetheless,  we obtain the required analysis using stochastic approximation based techniques (\cite{harold1997stochastic}): }

\Rev{\begin{thm}
\label{Thm_policy_smart_manu}
\textit{i) For any $\pi$ , $\br$ and $\bar D_{j,t}$ defined via \eqref{Eqn_dt_smart_manu}-\eqref{Eqn_q_phi_smart}  converges a.s. to the following unique limit:}

\vspace{-3mm}
\textit{{\small\begin{eqnarray}
     \bar d_j^*  \hspace{-2mm}&\hspace{-2mm}=\hspace{-2mm}&\hspace{-2mm}\frac{\sqrt{\alpha} \ \ \sum_{i=N_j(\bar d_j^*)}^L \sigma_i r_{ij}  \left(    h_{ij}     +  \pi \alpha h_{jj}     \right) }{ (\sum_{i=1}^L \sigma_i r_{ij} + \pi \alpha  \sum_{i=N_j(\bar d_j^*)}^L \sigma_i r_{ij}  )  }, \mbox{ \normalsize with, } \hspace{2mm}
     \label{Eqn_q_smart_manu}  \\
     N_j(\bar d) \hspace{-1.5mm}&\hspace{-1.5mm}:=\hspace{-1.5mm}&\hspace{-1.5mm} \min \left \{i:   h_{ij} > \sqrt{\alpha} \pi ({\bar d} -  \sqrt{\alpha}h_{jj}) \right  \} \ \mbox{({\small $\sum_{N_j}  = 0$} if empty set)}. \nonumber 
\end{eqnarray}}
Also the revenue rate of the manufacturer converges a.s.,}  

\vspace{-3mm}
\textit{
{\small
\begin{eqnarray}
{\bar U}_M^\infty (\br, \pi) &=& \lim_{t \to \infty} \frac{\sum_{s\le t}  D_s (P_s - Q_s - C_m) } {t} = \sum_{i,j} \sigma_i r_{ij} \bar d_j^*   W_{ij}^*,  \nonumber   \\
  W^*_{ij}  &:=&  {\sqrt{\alpha} h_{ij}- \pi \alpha ( \bar d^*_j   - \sqrt{\alpha} h_{jj} )  }. \label{Eqn_UM_conv}
\end{eqnarray}}
ii)   Under truthful reporting,  $\br = \brid$,  for any $\pi$, the \manu derives the `best' as in \eqref{Eqn_Um_star_limit_SBE}:}  

\vspace{-3mm}
\textit{
{\small \begin{eqnarray}
     {\bar U}_M^\infty (\brid, \pi):= {\bar U}_M^*   = \frac{1}{16 \alpha}\sum_i \sigma_i {\left([ \phi_i-\alpha (C_s + C_m)]\right)^2}. \mbox{ a.s.}
     \label{Eqn_rid_util} 
\end{eqnarray}} 
iii)    For any   $\br  \ne \brid$, there exists a $\bar \pi$ and  for any  $\pi > \bar \pi$,  we have,
$
 \bar U_M^\infty (\br, \pi) < \bar U_M^*.  
$}\eop
\end{thm}

Thus,  we establish that the policy \eqref{Eqn_q_phi_smart} based on the micro observations   works against the strategic manufacturer -- against any $\br \ne \brid$, when supplier uses   an appropriate $\pi$ in \eqref{Eqn_q_phi_smart} the strategic \manu derives strictly inferior; the \manu would derive the maximum as in \eqref{Eqn_Um_star_limit_SBE} if it instead reports truthfully (i.e., using $\brid$).  
\textit{The dynamic analysis for policy \eqref{Eqn_q_tilde} follows similarly as in Theorem \ref{Thm_policy_smart_manu}.}
}

\vspace{-2mm}
\section{Conclusions}
\label{sec_conclusion}
The agents of 
supply chains at various echelons have differential access to  important information, like instantaneous market potentials. This asymmetry can be utilized by the agents at the lower echelon for their advantage.  In this paper,  we study in depth a dynamic supply chain game subjected to the market fluctuations   with and without information asymmetry. We prove that the agent at the higher echelon can design punitive policies if it acts as the leader by first setting the prices. These policies can successfully extract truthful information from the greedy agents at the lower echelon that under-report the market potential. 
\Rev{We propose  a second punitive policy that  combats successfully (bur requires micro observations) a more strategic agent at the lower echelon that under and over reports to its advantage.}
This is just the beginning, one can investigate several future directions, analysis with estimation and prediction errors,   etc.  

\vspace{-2mm}
\section{Appendix}
\label{sec_app}
\vspace{-1mm}
\noindent \underline{\textbf{Proof of Theorem \ref{Thm_greedy_policy}:}}  
Define $\pi_{ij}:= \pi \sigma_i (\phi_i-\phi_j)$, $\psi_{ij} := \left [ h_{ij} -\pi f(\br)  \right ]^+$ and  ${\cal R}_{ij} (\pi) := \left \{\br \in \R_e : \pi f(\br)   \le h_{ij} \right  \}$, for all $i,j$. 
From \eqref{Eqn_h_ij},  $\Uij(i,j) = \psi_{ij}^2 > 0$, only for $\br $ in the interior of ${\cal R}_{ij} (\pi)$.
From \eqref{eqn_um_pi_rpi},   the partial derivative of $\tU_{M}(\br, \pi)$:

\vspace{-3mm}
{\small
    \begin{eqnarray}
  \frac{\partial \tU_{M}(\br, \pi)}{\partial r_{ij}} 
&=& \sigma_i\bigg( \Uij(i,j) - \Uij(i,i)\bigg) + \sum_{i'} \sigma_{i'} \frac{\partial \Uij(i',i')}{ \partial r_{ij}} \nonumber\\
 &&+ \sum_{(i' , j'): i' > j'} \sigma_{i'} r_{i' j} \left( \frac{\partial \Uij(i',j')}{ \partial r_{ij}} -\frac{\partial \Uij(i',i')}{ \partial r_{ij}} \right)  \nonumber 
 \\
  & & \hspace{-22mm} =\  \sigma_i  \left( \psi_{ij}^2 - \psi_{ii}^2 \right)  - 2  \pi_{ij} \left( \hspace{-1mm}\sum_{i'} \hspace{-1mm}\sigma_{i'}  \psi_{i' i'}  + \hspace{-3mm}\sum_{(i' , j'): i' > j' } \hspace{-5mm}\sigma_{i'} r_{i' j}  \left( \psi_{i' j'} - \psi_{i' i'} \right) \right). \hspace{7mm}
\label{Eqn_um_derivative_general}
\end{eqnarray}
}%
We prove the theorem using the following three cases: \\
{\bf Case 1:} If $\br \in  {\cal R}_{11}  $. In this  case $ \pi f(\br) \leq h_{11} < h_{LL}$ and so   the derivative w.r.t., $r_{ij}$ for any $i > j$ can be upper bounded:

\vspace{-4mm}
{\small    \begin{eqnarray*}
 \frac{\partial \tU_{M}(\br, \pi)}{\partial r_{ij}} &\le&   \sigma_i   \Uij(i,j) -  2 \pi_{ij}   \sigma_{L} \left(h_{LL} - \pi f(\br) \right ) \\
 &\le& \sigma_i   h_{ij}^2 -  2 \pi_{ij}   \sigma_{L} (h_{LL} -h_{11} ).
 \end{eqnarray*}}
Now, one can choose the $\pi$ (say $\pi_1$) large enough such that the above derivative remains strictly negative  for all $\br \in \R_{11}$ and for all $(i,j)$ (observe  $|\Q|<\infty$ and $h_{ij}$'s are all bounded); then $ \tU_{M}(\br, \pi) <  \tU_{M}(\brid, \pi)$ for $\br \in \R_{11}$. \\
{\bf Case 2:} Now consider the case when    $\br \in \R_{\ne}$, where 
\begin{eqnarray*}
\R_{\ne} \hspace{-1mm}&\hspace{-1mm}:=\hspace{-1mm}&\hspace{-1mm} \{\br \in \R_e : \Uij(i,i) = 0, \mbox{ for all } i\} = \R_e \setminus \R_{LL}. 
\end{eqnarray*} 
By extending the function suitably,  we have for any $\br \in \R_{\ne}$:

\vspace{-4mm}
{\small\begin{eqnarray}
    \tU_{M}(\br, \pi) \le  \max_{ \br \in \R_{\ne} } \sum_{ i > j} \sigma_i r_{ij} \Uij(i, j)  \stackrel{b}{\le}  \max_{ \br \in \R_e } \sum_{ i > j} \sigma_i r_{ij} \Uij(i, j), 
    \label{Eqn_R_hash}
 \end{eqnarray}}%
 where `$b$' is true since $\R_{\ne}  \subset \R_{e} $. 
 The maximum in the last term   is attained at a $\br^*$ with
   $ \pi f(\br^*)  \le  h_{L1}$ 
 (as $h_{L1} \ge h_{ij}$   $\forall (ij)$, so  the objective function in \eqref{Eqn_R_hash} is $0$ $\forall \br$ with $ \pi f(\br)  \ge  h_{L1}$). 
In other words, the  components $r_{ij}^*$ of $\br^*$ satisfy

\vspace{-2mm}
{\small$$
r_{ij}^* \leq \frac{h_{L1}}{\pi \sigma_i (\phi_i - \phi_j)}, \mbox{ for all } i >  j.
$$}%
So from
 \eqref{Eqn_R_hash},  $
      \tU_{M}(\br, \pi) \le \frac{1}{\pi} \sum_{i>j}\frac{h_{L1}}{(\phi_i - \phi_j)} (h_{ij})^2, 
 $ for all $\br \in \R_{\ne} $ (note $\Uij(i,j)\le h_{ij}^2$ $\forall (ij)$).
 Thus   choosing sufficiently large $\pi$ (say $\pi_2$) ensures  $
      \tU_{M}(\br, \pi) <  \tU_{M}(\brid, \pi)$ for $\br \in \R_{\ne}$.  \\
{\bf Case 3:} Now consider the left-over region, $\br \in \R_{11}^c \setminus \R_{\ne}$.
For such $\br$, we have $\Uij(1,1)=0$ and hence, 

\vspace{-3mm}
{\scriptsize \begin{eqnarray*}
    \tU_{M}(\br, \pi) &=& \sum_{i>j} \sigma_i r_{ij} \Uij(i,j) + \sum_{l \ge 2} \sigma_l \Uij(l,l) 
    \\
    & & \hspace{-20mm}\leq \  \max_{ \br \in \R_e } \sum_{ i > j} \sigma_i r_{ij} \psi_{ij}^2 + \sum_{l \ge 2} \sigma_l \Uij (l,l; \brid) \ <  \  \frac{ \sum_{i>j}\frac{h_{L1}}{(\phi_i -\phi_j)} (h_{ij})^2 }{\pi}+ \sum_{l \ge 2} \sigma_i \Uij(l,l;\brid). 
\end{eqnarray*}}%
Using logic as in Case 2,   choose large $\pi$ (say $\pi_3$), and then 

\vspace{-3mm}
{\small\begin{eqnarray*}
  \frac{1}{\pi}  \sum_{i>j}\frac{h_{L1} (h_{ij})^2}{(\phi_i - \phi_j)}  + \sum_{l \ge 2} \sigma_i \Uij (l,l;\brid) \  < \sum_{l \ge 1} \sigma_i \Uij (l,l;\brid) = U_M(\brid, \pi). 
\end{eqnarray*}}
Now, consider $\bar{\pi} = \max\{\pi_1, \pi_2, \pi_3 \}$ and hence the proof. \eop

\TR{
\Rev{\noindent \underline{\textbf{Proof of Theorem \ref{Thm_policy_smart_manu}, Part (i):}} 
 Consider  $\phi_j$  such that   $P(\hPhi_1 = \phi_j)>0$. Since  sequence $\{\hPhi_t \}_t$ is IID,   to derive asymptotic analysis of  $\bar D_{j,t}$,  it suffices  to consider  subsequence $\{t_k\}$ such that $\hPhi_{t_k} = \phi_j$. For simplicity, rename $\{t_k\}$ by  $\{t\}$ and then $\bar D_{j,t}$ of \eqref{Eqn_report_demand_smart_manu} evolves as (see $D_t$ in \eqref{Eqn_dt_smart_manu}),
 
 \vspace{-3.5mm}
 {\small
 \begin{eqnarray*}
 \hspace{1mm}
    \bar D_{j, t} = \bar D_{j, (t-1)} 
+ \nicefrac{1 }{t} \left(D_t - \bar D_{j, (t-1)} \right), \ \  \mbox{\normalsize along the subsequence. }  
\end{eqnarray*}}%
Using stochastic approximation result of   \cite{harold1997stochastic}, we would   approximate  $\{\bar D_{j,t}\}$   by the following  ODE (see \eqref{Eqn_pstar_qstar_limit_SBE}, \eqref{Eqn_q_phi_smart}): 

\vspace{-2mm}
{\small\begin{eqnarray}
\label{Eqn_ODE_for_Dj}
 \dot {\bar d}_j  =  {\bar g} (\bar d_j)   &:= & E[D_t- \bar{D}_{j, (t-1)}|   \bar D_{j,  (t-1)} = \bar d_j, \hPhi_t = \phi_j]    \\
  & & \hspace{-25mm} =\   \frac{  
  \sum_{i = N_j(\bar d_j)}^L \sigma_i r_{ij}  (\sqrt{\alpha }h_{ij}- \pi \alpha   \left( {\bar { d_j}}  -  \sqrt{\alpha }h_{jj} \right) ) }{   \sum_{i=1}^L \sigma_i r_{ij}}    - \bar d_j, \  (N_j\mbox{ as in  \eqref{Eqn_q_smart_manu}} ).  \nonumber
\end{eqnarray}}%
Towards proving the above, first observe ${\bar g}  $  is  strictly decreasing and Lipschitz continuous, has unique zero which is the unique attractor  of~\eqref{Eqn_ODE_for_Dj}:
\begin{eqnarray}
      \label{Eqn_global_attractor}
     \bar d_j^* =   \bar d_j^*(\br; \pi) =  \frac{\sqrt{\alpha}\sum_{i=N_j(\bar d^*_j)}^L \sigma_i r_{ij}  ( h_{ij}+ \alpha \pi h_{jj} ) }{    \sum_{i=1}^L \sigma_i r_{ij} +  {\pi \alpha \sum_{i=N_j(\bar d^*_j)}^L \sigma_i r_{ij} } }; 
\end{eqnarray}
note  $\bar g(0) >0$ and $\bar g(\phi_L) < 0$ and hence
the above follows by intermediate value theorem, see \cite{arxiv} for more details.
Further,  it is easy to verify assumptions A2.1-A2.5
of \cite{harold1997stochastic}.
%
  %
Thus part (i) follows by \cite[Theorem 2.1, pg. 127]{harold1997stochastic} (clearly $P (\bar D_{j,t} \in [0, \phi_L] \mbox{ i.o.}) = 1 $ and almost sure convergence of  \eqref{Eqn_UM_conv} follows from that of $\bar D_{j,t}$ (for each $j$) using standard arguments,   again \cite{arxiv} has more details).  

\noindent\underline{\bf Part (ii):}
When $\br = \brid$, from \eqref{Eqn_global_attractor},  the unique 
attractor  
 equals (observe $r_{ij} = 0$ when $i\ne j$ here)
$ 
\bar d_j^*(\brid; \pi) = \sqrt{\alpha} h_{jj} 
$ for any $j$. So from \eqref{Eqn_UM_conv}, the asymptotic utility equals \eqref{Eqn_rid_util}.


\noindent\underline{\bf Part (iii):} 
From \eqref{Eqn_global_attractor},  we have 
  $\bar d_j^*(\br; \pi) \to \sqrt{\alpha} h_{jj}$ for any $\br, j$, as $\pi \to \infty$, and so
{\small$ \bar V  (\br, \pi) := \sum_{i,j} \sigma_i r_{ij} \bar d_j^*(\br; \pi) \sqrt{\alpha} h_{ij} \to \alpha \sum_{i, j} \sigma_i r_{ij}  h_{jj} h_{ij} . $}  Further from \eqref{Eqn_UM_conv},   {\small $ \bar U_M^\infty (\br, \pi) < \bar V  (\br, \pi) $} and

\vspace{-3mm}
{\small
\begin{eqnarray*}
 \lim_{\pi \to \infty}  \bar V  (\br, \pi) -  \bar U_M^\infty (\brid, \pi)  
& = & \alpha\sum_{i,j} \sigma_i r_{ij}    h_{jj} h_{ij}
- \alpha \sum_i \sigma_i  h_{ii}^2  \sum_{j} r_{ij} \\
&= & \frac{ \sum_i \sigma_i \sum_j r_{ij}  \left(-( \phi_j - \phi_i)^2  \right)  }{16}  < 0. \hspace{1mm} \mbox{\eop}
\end{eqnarray*}}

\vspace{-2mm}


}}{}

\vspace{-2mm}
\bibliographystyle{IEEEtran}

\section*{Appendix for  Stochastic Approximation}
\label{App_for_sa}

\underline{\bf Proof of Theorem \ref{Thm_policy_smart_manu}, Part (i):}
Observe that $\bar D_{j, t}$ gets updated only when 
$\hPhi_t = \phi_j$,  and also observe under SMR policies $(\br, \pi)$, 
 $\{\hPhi_t\}_t$ is an IID sequence. Thus,  for any $\phi_j$  we have:
\begin{eqnarray}
    P(\hPhi_t = \phi_j) =  P(\hPhi_1 = \phi_j), \mbox{ for all }  t. 
    \label{Eqn_pob_phi_positive}
\end{eqnarray}
Define $
N_{jt} := \sum_{s=1}^t \mathds{1}_{\{\hPhi_s = \phi_j \}}$ for any  $\phi_j$, then 
by law of large numbers, $  \nicefrac{N_{jt}}{t} \to P(\hPhi_1 = \phi_j ) $ almost surely. Now consider  a  $\phi_j$ with $P(\hPhi_1 = \phi_j ) > 0$  and then we have:
$$
N_{jt} = \sum_{s=1}^t \mathds{1}_{\{\hPhi_s = \phi_j \}} \to \infty,  \mbox{ almost surely as }  t\to \infty.
$$

We now consider  such $\phi_j$ and focus only on the subsequence $\{N_{j,k_t}\}_t$  at which    $\phi_j$ is reported, to obtain time-asymptotic analysis of $\bar D_{j, t}$ (it is clear that $\bar D_{j', t}$ is not relevant for analysis if  $P(\hPhi_1 = \phi_{j'}) = 0$).  \textit{For simpler notation, we represent subsequence $\{k_t\}$ by $\{t\}$ and focus only on the events pertaining to those epochs.} By virtue of the  simpler notation,   we have $\hPhi_t = \phi_j$, for all $t$. 

Define filtration $\mathcal{F}_t^{j}:= \sigma(D_s: s \le t)$, which  represents the sigma-algebra generated by $(D_s: s \le t)$, where $D_s$ is given by \eqref{Eqn_dt_smart_manu} and \eqref{Eqn_q_phi_smart}. 
One can write $\bar D_{jt}$ in  recursive manner as below (for the sub-sequence under focus):
\begin{eqnarray}
    \bar D_{j, t} = \bar D_{j, (t-1)} 
+ \frac{1 }{t} (D_t - \bar D_{j, (t-1)} ). 
\label{Eqn_stoc_Approx}
\end{eqnarray}
This represents a stochastic approximation based iterative scheme (e.g., as in \cite{harold1997stochastic})  depending upon the the policy \eqref{Eqn_q_phi_smart} and the limit can be approximated by the following ODE   (see also  \eqref{Eqn_pstar_qstar_limit_SBE}, \eqref{Eqn_h_ij})
\begin{eqnarray}
\dot {\bar d}_j  &=&  {\bar g} (\bar d_j)  \mbox{, where, }\label{Eqn_ODE_for_Dj_new}\\
{\bar g} (\bar d)  &:=& \ E[D_t- \bar{D}_{j, (t-1)}| \mathcal{F}_t^{j}, \bar D_{j,  (t-1), } = \bar d, \hPhi_t = \phi_j] \nonumber \\ 
&  & \hspace{-20mm} = \ \frac{  
  \sum_{i = N_j(\bar d)}^L \sigma_i r_{ij}  (\sqrt{\alpha }h_{ij}- \pi \alpha   \left( {\bar { d}}  -  \sqrt{\alpha }h_{jj} \right) ) }{   \sum_i \sigma_i r_{ij}}    - \bar d,  \label{Eqn_bar_g_bar_d}  \\
\mbox{ with, } && N_j(\bar d) := \min \left \{i:   h_{ij} > \sqrt{\alpha} \pi ({\bar d} -  \sqrt{\alpha}h_{jj}) \right  \}.  \nonumber 
\end{eqnarray}
It is not difficult to show that the scheme in \eqref{Eqn_stoc_Approx}
satisfies the assumptions A2.1-A2.5  of \cite[Theorem 2.1, pg. 127]{harold1997stochastic}
with $\beta_n, \epsilon_n$ and $\bar g$ of \cite{harold1997stochastic} in our case given by:
$$
\beta_n = 0, \ \epsilon_n = \frac{1}{n} \mbox { and } \  {\bar g}  \mbox{ of \eqref{Eqn_bar_g_bar_d}, respectively} . 
$$
From  \eqref{Eqn_h_ij},   the set $[0,   \phi_L]$ is invariant  for the ODE \eqref{Eqn_ODE_for_Dj_new} as $\bar g(\phi_L  ) < 0$, recall $\phi_i \le \phi_L$ for all $i$. 
From  \eqref{Eqn_dt_smart_manu}-\eqref{Eqn_report_demand_smart_manu} one can upper bound (recall $\{{\cal N}_s\}$ bounded IID zero mean noise), 
\begin{eqnarray*}
    \bar D_{j, t}  \le \frac{ \sum_{s \le t}  \Phi_s + {\cal N}_s } {t} 
\end{eqnarray*}
and the upper bound converges  by  LLNs, almost surely,
$$
\frac{ \sum_{s \le t}  \Phi_s + {\cal N}_s } {t}  \to  \phi_L \mbox{ almost surely as }  t \to \infty. 
$$
Hence,  $P( \bar D_{j, t} \in [0,  \phi_L ] \mbox{ i.o.}) = 1$,  for the scheme in \eqref{Eqn_stoc_Approx}.  

We will now show that the ODE \eqref{Eqn_ODE_for_Dj_new}  has unique global attractor and then the proof of  part (i) is completed by
   \cite[Theorem 2.1, pg. 127]{harold1997stochastic}. 
    
{\bf Proof of existence of unique global attractor:}
To begin with RHS of the ODE,   $\bar g$ is Lipschitz continuous  
and  linear, which proves the existence of unique ODE solution for any given initial condition. Towards further proof, observe the RHS of ODE from \eqref{Eqn_bar_g_bar_d}: (i) the  second term $(-\bar d_j)$ is   decreasing in $\bar d_j$;    (ii) mapping $\bar d \mapsto N_j(\bar d)$ is piecewise constant and  non-decreasing,   eventually the set defining $N_j(\bar d)$ becomes empty, implying the first term becomes zero;  and finally (iii) each individual term in the summation of  the first term in RHS is also decreasing with~$\bar d_j$.   

Thus,  RHS $\bar g$  is a strictly decreasing function of $\bar d$,   $\bar g(0) >0 $ and $\bar g(\phi_L)<0$.  
Thus, further using  intermediate value theorem, there exists a unique zero  in $[0, \phi_L]$ which becomes the global  attractor and from \eqref{Eqn_bar_g_bar_d}, it is given by:
\begin{eqnarray}
    \label{Eqn_global_attractor_new}
\bar d_j^*(\br; \pi) =  \frac{\sqrt{\alpha}\sum_{i=N_j(\bar d_j^*)}^L \sigma_i r_{ij}  ( h_{ij}+ \alpha \pi h_{jj} ) }{   \sum_i \sigma_i r_{ij} +  {\pi \alpha \sum_{i=N_j(\bar d_j^*)}^L \sigma_i r_{ij} } }. 
\end{eqnarray}
As this is global attractor,  $[0, \phi_L]$ is in the (domain of attraction) DoA corresponding to $\bar d_j^*(\br; \pi)$.

{\bf Almost sure  convergence in \eqref{Eqn_UM_conv}:}  Using the definition of  $Q_s $   in \eqref{Eqn_q_phi_smart},  we have for any time slot $s$: 
$$
Q_s  =    q^* (\hPhi_s)  + 2 \pi ( \bar D_{\hPhi_s, s-1} - \sqrt{\alpha} h_{jj} ).  
$$
Thus by  almost sure convergence  of $\bar D_{j, s-1} \to \bar d^*_j (\br, \pi)$  (for each~$j$) we have:
 \begin{eqnarray}
 \left | Q_s- Q_s^* \right  | & \stackrel{s\to \infty}{\to} &  0  \mbox{ a.s., where }\label{Eqn_Qs} \\
 Q_s^* & := & \hspace{-4mm} \left ( \sum_{j} \mathds{1}_{\{\hPhi_s = \phi_j\}} q^* (\phi_j)  + 2 \pi ( \bar d^*_j (\br, \pi) - \sqrt{\alpha} h_{jj} )  \right ). \nonumber 
 \end{eqnarray} 
 Recall $P_s   =  \nicefrac{(\Phi_s+\alpha(Q_s+C_m))}{2\alpha}$ for any $s$ and  similar convergence can be derived for $P_s$. Consider the sample paths in which $P_s$ and $Q_s$ have converged (this happens with probability one). Fix one such sample path,  there exists a $\bar t_\epsilon$ such that 
 \begin{eqnarray}
 \left | Q_s-  Q^*_s \right |  
   \le  \epsilon \mbox{ and }  \left | P_s-  P^*_s \right |  
   \ \le  \ \epsilon, \mbox{ for all } s \ge \bar t_\epsilon, 
   \label{Eqn_p_q_le_epsi}
   \end{eqnarray}
    where (represent $\bar d_j^* (\br, \pi)$   as $\bar d_j^*$ for brevity, when there is clarity) 

\vspace{-2mm}
    {\small
  \begin{eqnarray}
  P^*_s &=&  \hspace{-2mm}
   \sum_{i,j}  \mathds{1}_{\{\hPhi_s = \phi_j, \Phi_s = \phi_i \}} \frac{\phi_i + \alpha( q^*(\phi_j) + 2 \pi ( \bar d^*_j  - \sqrt{\alpha} h_{jj} ) +C_m)}  {2\alpha}  \nonumber \\
  &=&  \hspace{-2mm} 
   \sum_{j} \mathds{1}_{\{\hPhi_s = \phi_j, \Phi_s = \phi_i\}} \frac{2 \phi_i + \phi_j + \alpha (C_s + C_m)+ 2 \pi \alpha ( \bar d^*_j  - \sqrt{\alpha} h_{jj} ) }{4 \alpha }    \nonumber \\
   \label{Eqn_Ps}
 \end{eqnarray}}
 Define  
 \begin{eqnarray}
      W^*_{ij}  := {\sqrt{\alpha} h_{ij}- \pi \alpha ( \bar d^*_j  (\br, \pi) - \sqrt{\alpha} h_{jj} )  }, \label{Eqn_Wij_star}
 \end{eqnarray}
 and observe from \eqref{Eqn_Ps} and \eqref{Eqn_Qs} we have:
 \begin{eqnarray}
 \label{Eqn_PsQs_star}
    P_s^* - Q_s^* - C_m = \sum_{i,j}  \mathds{1}_{\{\hPhi_s = \phi_j, \Phi_s = \phi_i\}} W^*_{ij}.
 \end{eqnarray}

Now one can write the terms defining limit \eqref{Eqn_UM_conv} as:
 \begin{eqnarray}
       \frac{\sum_{s\le t}  D_s (P_s - Q_s - C_m) } {t} 
       = \frac{\sum_{s\le t}  D_s ( (P_s - Q_s - P^*_s + Q_s^*)} {t} \nonumber \\
       +  \frac{\sum_{s\le t}  D_s (  P^*_s - Q^*_s - C_m ) } {t} . \label{Eqn_Ds_conv}
       \end{eqnarray}
The first term in the above converges to zero, as:
\begin{itemize}
    \item 
for all $t \ge \bar t_\epsilon$, using \eqref{Eqn_p_q_le_epsi} (observe all terms are upper bounded a.s. by a uniform bound, say $B_1$), 

\vspace{-2mm}
{\scriptsize
 \begin{eqnarray*}
  \left |  \frac{\sum_{s\le t}  D_s ( P_s - Q_s -  P^*_s + Q^*_s  ) } {t}  \right | \hspace{-20mm} \\
  \hspace{-1mm}&\hspace{-1mm}=\hspace{-1mm}&  \hspace{-1mm} \left |\frac{\sum_{s<  \bar t_\epsilon}  D_s ( P_s - Q_s -  P^*_s + Q^*_s  ) + \sum_{s=   \bar t_\epsilon}^t  D_s ( P_s - Q_s -  P^*_s + Q^*_s  )}{t} \right | \\
  \hspace{-2mm}&\hspace{-2mm}\le \hspace{-2mm} & \hspace{-2mm} \left |\frac{\sum_{s<  \bar t_\epsilon}  D_s ( P_s - Q_s -  P^*_s + Q^*_s  )}{t} \right | + \left |\frac{\sum_{s =  \bar t_\epsilon}^t  D_s ( P_s - Q_s -  P^*_s + Q^*_s  )}{t} \right | \\
  \hspace{-2mm}&\hspace{-2mm}\le \hspace{-2mm} & \hspace{-2mm} \frac{ B_1 \  \bar t_\epsilon} {t} +  2 \left |   \frac{\sum_{s =   \bar t_\epsilon}^t  D_s } {t}  \right | \epsilon,    \\
  \hspace{-2mm}&\hspace{-2mm}\le \hspace{-2mm} & \hspace{-2mm} \frac{ B_1 \  \bar t_\epsilon} {t} +  2 B_2 \epsilon, 
  \end{eqnarray*}
  }%
 the last inequality is because $D_s$ is a.s. upper bounded  by a positive constant, say  $B_2$;
 
 \item thus, we have 
 $$
0 \le \lim \sup_{t \to \infty}  \left |  \frac{\sum_{s\le t}  D_s ( P_s - Q_s - P^*_s + Q^*_s ) } {t}  \right |  \le 2 B_2  \epsilon;
 $$ 

 \item finally, the above is true for any $\epsilon > 0$, and this implies that  the limit exists and  equals zero:
 $$
  \lim_{t \to \infty}  \left |  \frac{\sum_{s\le t}  D_s ( P_s - Q_s - P^*_s + Q^*_s ) } {t}  \right | = 0.$$
\end{itemize}

Now consider the second term in \eqref{Eqn_Ds_conv} and using \eqref{Eqn_PsQs_star}  we have the following convergence (see~\eqref{Eqn_global_attractor_new})

\vspace{-2mm}
{\small
 \begin{eqnarray}
     \frac{\sum_{s\le t}  D_s (  P^*_s - Q^*_s - C_m ) } {t} &=& \hspace{-3mm}  \frac{\sum_{i,j} \sum_{s\le t}   \mathds{1}_{\{\hPhi_s = \phi_j, \Phi_s = \phi_i \}} D_s  W_{ij}^* } {t} \nonumber\\ 
     & & \hspace{-7mm} \stackrel{t\to \infty}{\rightarrow} \sum_{i,j} \sigma_i r_{ij} \bar d_j^* W_{ij}^*,    \hspace{2mm}
     \label{Eqn_second_term_limit}
 \end{eqnarray}}
and the above is true as for any $i,j$ 
 \begin{eqnarray*}
\frac{  \sum_{s\le t}   \mathds{1}_{\{\hPhi_s = \phi_j, \Phi_s = \phi_i\}} D_s W_{ij}^* }{t} \hspace{-30mm} \\
&=&  W_{ij}^* \frac{\sum_{s\le t}   \mathds{1}_{\{\hPhi_s = \phi_j, \Phi_s = \phi_i\}}}{t}
\frac{  \sum_{s\le t}   \mathds{1}_{\{\hPhi_s = \phi_j, \Phi_s = \phi_i\}} D_s  }{\sum_{s\le t}   \mathds{1}_{\{\hPhi_s = \phi_j, \Phi_s = \phi_i\}}} \\
&\stackrel{t\to \infty}{\rightarrow} & 
 W_{ij}^* \ \sigma_i r_{ij} \  \bar d_j^*;   
 \end{eqnarray*}
observe $\{ \mathds{1}_{\{\hPhi_s = \phi_j, \Phi_s = \phi_i\}} \}_s$ is a further subsequence in which $\phi_j$ is reported and actual potential $\Phi_s = \phi_i$ and (almost sure) convergence among main sequence implies the same along this subsequence:
 
 Thus, from \eqref{Eqn_Ds_conv}-\eqref{Eqn_second_term_limit}, we have  (see \eqref{Eqn_Qs})
 \begin{eqnarray}
 \frac{\sum_{s\le t}  D_s (  P_s - Q_s - C_m ) } {t} \stackrel{t\to \infty}{\rightarrow}  \sum_{i, j} \sigma_i r_{ij} \bar d_j^* W_{ij}^* . \hspace{4mm} 
     \label{Eqn_Ds_limit_final}
 \end{eqnarray}

\underline{\bf Proof of part (ii):} 
When $\br = \brid$, {\small$P(\hPhi_1 =\phi_j) = \sigma_j >0$}   $\forall  j$. From \eqref{Eqn_global_attractor_new},  the unique 
attractor  
 equals (observe $r_{ij} = 0$ when $i\ne j$ here):
$$  
\bar d_j^*(\brid; \pi) = \sqrt{\alpha} h_{jj} .
$$ 
And this is true for each $j$ and hence the asymptotic utility equals \eqref{Eqn_rid_util}, see \eqref{Eqn_Ds_limit_final} and \eqref{Eqn_Wij_star}. 

\underline{\bf Proof of  part (iii):}
From \eqref{Eqn_global_attractor_new}, for all $\br$ and $j$, we have:
\begin{eqnarray*}
    \bar d_j^*(\br; \pi)  &=&  \frac{\sqrt{\alpha}\sum_{i=N_j(\bar d_j^*)}^L \sigma_i r_{ij}  \left( \frac{h_{ij}}{\pi}+ \alpha h_{jj} \right) }{  \frac{ \sum_i \sigma_i r_{ij}}{\pi} +   \alpha \sum_{i=N_j(\bar d_j^*)}^L \sigma_i r_{ij} }  \\
     & \to & \sqrt{\alpha} h_{jj}, \mbox{ as } \pi \to \infty. 
\end{eqnarray*}
Define
$$
\bar V  (\br, \pi) := \sum_{i,j} \sigma_i r_{ij} \bar d_j^* \sqrt{\alpha} h_{ij}
$$
and observe from \eqref{Eqn_Wij_star} that $ \bar U_M^\infty (\br, \pi) < \bar V  (\br, \pi) $ for all $\br \ne \brid$. 
Now from \eqref{Eqn_Ds_limit_final}, as $\pi \to \infty$ we have:
$$  \bar V  (\br, \pi)  \to  \alpha \sum_{i, j} \sigma_i r_{ij}  h_{jj} h_{ij} . $$ 
Now, 
\begin{eqnarray*}
 \lim_{\pi \to \infty}  \bar V  (\br, \pi) -  \bar U_M^\infty (\brid, \pi)   
& = & \alpha\sum_{ij} \sigma_i r_{ij}    h_{jj} h_{ij}
- \alpha \sum_i \sigma_i  h_{ii}^2  \sum_{j} r_{ij} \\
&=& \alpha \sum_i \sigma_i \sum_j r_{ij} \left (h_{jj}  h_{ij} - h_{ii}^2 \right ) \\ 
&\stackrel{b}{=}& \frac{ \sum_i \sigma_i \sum_j r_{ij}  \left(-( \phi_j - \phi_i)^2   \right)  }{16}  < 0,
\end{eqnarray*}%
where equality `b' follows as:
\begin{eqnarray*}
   16 \alpha( h_{jj}  h_{ij} - h_{ii}^2) &=& (\phi_j - \alpha(C_s
    + C_m)) (2 \phi_i - \phi_j - \alpha(C_s
    + C_m) ) \\
    && - (\phi_j - \alpha(C_s
    + C_m))^2 \\
    &=& 2 \phi_i \phi_j - \phi_j^2 -  \alpha \phi_j (C_s
    + C_m) -  2\alpha \phi_i (C_s+ C_m) \\
    &&+ \alpha \phi_j (C_s
    + C_m) + \alpha^2 (C_s+ C_m)^2 \\
    &&- \left( \phi_i^2 + \alpha^2 (C_s+ C_m)^2- 2 \alpha \phi_i (C_s+ C_m) \right) \\
    &=& 2 \phi_i \phi_j - \phi_j^2 - \phi_i^2 = - (\phi_j- \phi_i)^2 
\end{eqnarray*}
Hence the result.  \eop

\end{document}